\documentclass[11pt]{article}

\title {On the Combinatorial Structure of Primitive\\Vassiliev Invariants, III -
A Lower Bound}
\author {Oliver T. Dasbach 
\thanks{supported by the Deutsche Forschungsgemeinschaft (DFG)}
\thanks{e-mail: \sl kasten@math.columbia.edu, 
http://www.math.uni-duesseldorf.de/home/kasten}\\\
{\small \sl Columbia University}\\
          {\small \sl Department of Mathematics}\\
          {\small \sl New York, NY 10027}
	  }

\date{}
\input{epsf}
\newtheorem{satz}{Satz}[section]
\newtheorem{example}[satz]{Example}
\newtheorem{lemma}[satz]{Lemma}
\newtheorem{definition}[satz]{Definition}
\newtheorem{proposition}[satz]{Proposition}
\newtheorem{theorem}[satz]{Theorem}
\newtheorem{corollary}[satz]{Corollary}

\newtheorem{maintheorem}[satz]{MAIN THEOREM}

\newenvironment
  {conjecture} {\bigskip \noindent{\bf Conjecture }\it}{\bigskip}

\newenvironment
  {proof}{\noindent {\bf Proof }}{$\Box$\bigskip \medskip}

\newenvironment
  {remark}{\bigskip \noindent {\bf Remark }}{\bigskip}

\newcommand {\B} {{\mathscr{B}}}
\newcommand {\W} {{\mathscr W}}

\newcommand {\A} {{\mathscr {A}}}

\newcommand {\V} {{\mathscr V}}

\newcommand {\IHXr} {{{\it IHX}-relation}}
\newcommand {\asr} {{antisymmetry relation}}

\newcommand {\C} {\mathbb{C}}

\newcommand{\ihx}
{\unitlength=0.5mm
\begin{picture}(62.00,13.00)
\put(0.00,13.00){\line(1,0){16.00}}
\put(0.00,1.00){\line(1,0){16.00}}
\put(7.00,13.00){\line(0,-1){12.00}}
\put(19.00,7.00){\makebox(0,0)[cc]{$=$}}
\put(26.00,13.00){\line(0,-1){12.00}}
\put(26.00,6.00){\line(1,0){12.00}}
\put(38.00,13.00){\line(0,-1){12.00}}
\put(43.00,7.00){\makebox(0,0)[cc]{$-$}}
\put(53.00,4.00){\line(1,0){6.00}}
\put(50.00,1.00){\line(1,1){12.00}}
\put(62.00,1.00){\line(-1,1){12.00}}
\end{picture}
}

\usepackage{latexsym}
\usepackage{mathrsfs}
\usepackage{amssymb}

\newcommand {\Lie}{{\mathscr L}}

\newcommand {\gl} {{\mathfrak{gl}}}

\newcommand {\cP} {{\mathscr P}}

\newcommand {\oWgl} {{\overline{\W}_{\!\gl(N)}}}

\begin{document}


\maketitle

\begin{abstract}
We prove that the dimension of the space of primitive
Vassiliev invariants of degree $n$ grows - as $n$ tends to infinity -
faster than $e^{c \sqrt n}$ for any $c< \pi \sqrt {2/3}$.

The proof relies on the use of the weight systems coming from the Lie
algebra $\gl(N)$. In fact, we show that our bound is - up to a multiplication 
with
a rational function in $n$ - the best possible that one can get with
$\gl(N)$-weight systems. 
\end{abstract}

\section{Introduction}

The space $\V$ of Vassiliev knot invariants is still mysterious.
Although we have a perfect combinatorial description of it
(see e.g. \cite{BL,Kontsevich,BarNatan1}),
even the asymptotic behavior of the dimension of Vassiliev invariants
in degree $n$ is unknown.

Vassiliev invariants form an algebra which
is isomorphic to a free polynomial algebra.  A basis for this algebra is given by
a basis of the primitive Vassiliev invariants.
Therefore all information about Vassiliev invariants is included in
the primitive ones.
Soon after the discovery of Vassiliev invariants, it became clear that there is at 
least one
primitive Vassiliev invariant in each degree \cite{CDL3}. From that it
follows that the dimension of the space of (not necessarily primitive) Vassiliev 
invariants of order
$n$ grows - as $n$ tends to infinity -
faster than $e^{c \sqrt n}$ for any $c< \pi \sqrt {2/3}$
\cite{Kontsevich}.

Other subspaces of the primitive space can be obtained in the following
way:

For an arbitrary subspace of Vassiliev invariants one can take the
subalgebra generated by this subspace and look at the intersection in
degree $n$ with the primitive space.
One can see \cite{CDV} that this construction gives for
the HOMFLY-Vassiliev invariants \cite{BL} a contribution to the primitive
space of dimension $\left \lfloor \frac n 2 \right \rfloor$ in degree $n$.
The dimension of the subspace of the primitive space coming in this way
from the (unframed) colored Jones polynomial \cite{MM} in degree $n$ is again
$\left \lfloor \frac n 2 \right \rfloor$ \cite{CDV}.

The best known lower bound for the dimension of the primitive
space was recently given in a nice paper by Chmutov and Duzhin \cite{CD2}.
They proved that the primitive space in degree $n$ is at least
$n^{\log n}$ dimensional as $n$ tends to infinity.

The first aim of this paper is to prove the so called
Kontsevich-Bar-Natan conjecture (for the history of it see \cite{CD2})
which states that
the dimension of the space of primitive
Vassiliev invariants of degree $n$ grows - as $n$ tends to infinity -
faster than $e^{c \sqrt n}$ for any $c< \pi \sqrt {2/3}$.
Therefore we will get a much better lower bound than $n^{\log n}$.

For showing this we will make use of the universal weight system coming from
the Lie algebra $\gl(N)$, which is related to the cablings of the
HOMFLY-Vassiliev invariants.

As a by-product of the proof we will get that
the vector space $C_a H(n)$ of all Vassiliev invariants that come
(in the usual way) from the HOMFLY-Vassiliev invariants in degree $n$
and all of their cablings, connected and disconnected,
behaves like
$$g_1(n) e^{c \sqrt n} \leq \dim C_a H(n) \leq g_2(n) e^{c \sqrt n}$$
for two rational
functions $g_1$ and $g_2$ in $n$ and $c=\pi \sqrt {2/3}$.

I would like to express my deep gratitude to Sergei Chmutov and Arkady Vaintrob
for some months of continuous and fruitful discussions during their visits
at the Max-Planck-Institut, Bonn.

Furthermore I both thank Joan Birman as well as Gregor Masbaum for useful remarks on
an earlier draft.

This paper covers a talk given at the "Knot theory week", Bonn, July 1997,
well-organized by C.-F. B\"odigheimer and his knot theory group.

In a fax to Sergei Chmutov \cite{KontsevichFax}, Maxim Kontsevich 
independently gave a proof
for a weaker lower bound estimate. His bound is - roughly speaking - the square root
of ours. He uses the same techniques that we use.

\section{Preliminaries}

Vassiliev invariants form a filtered vector space $\V = \bigcup \V_n$ with
$\V_n$ the space of Vassiliev invariants of order at most $n$.
Kontsevich (see \cite{Kontsevich,BarNatan1}) gave a combinatorial
description of $\V_n/\V_{n-1}$:

The algebra $\A$ is the algebra generated by all chord-diagrams
modulo the four-term-relation,  with the connected sum as a multiplication.
It is graded by the number of chords in a diagram and the graded n-part
of it is denoted by $\A_n$.
By the result of Kontsevich we know that $\V_n/\V_{n-1}$ is isomorphic
to the subspace of the dual $\A_n^*$ of $\A_n$ that is generated by all 
functionals
vanishing on chord diagrams with an isolated chord.
Elements of $\A_n^*$ are called weight systems.

It turns out that with a suitable coproduct $\A$ becomes
an associative, commutative, coassociative and cocommutative Hopf-algebra.
By the classical structure theory of these algebras we know that $\A$ is
isomorphic to the polynomial algebra over its primitive space $\cP(\A)$.
The space $\cP(\A)$ corresponds to the subspace of $\V$ generated by invariants
$v$ that are additive for connected sums of knots:
$v(K_1 \# K_2)=v(K_1)+v(K_2)$.

We have another description of $\A$:
The algebra $\B$ is the algebra of all (finite) diagrams (graphs) having only
trivalent and univalent vertices, each trivalent vertex equipped with
one of the two cyclic orientations.
Furthermore the following two types of relations hold:
\begin{enumerate}
\item The \IHXr : \parbox{4cm}{$\ihx$}
\item The \asr : if in a diagram $D$ the orientation at one trivalent
vertex is changed then $D$ changes to $-D$.
\end{enumerate}

As usual, in all pictures of diagrams (or subdiagrams) in $\B$, it
is assumed that the three edges meeting at one trivalent vertex are
oriented counterclockwise.

The gradation in $\B$ is given by half of the number of vertices in a
diagram. The $n$-graded part is denoted by $\B_n$.
As vector spaces $\A_n$ and $\B_n$ are isomorphic and we have two 
different products
in $\B$: the natural product induced by the disjoint union of
diagrams and the product coming from the product in $\A$.

Since the \IHXr\ and the \asr\ are homogenous in the number of univalent
vertices we get a splitting $\B = \bigoplus \B^{(u)}$ where
$\B^{(u)}$ is the subspace generated by all diagrams in $\B$ with $u$ univalent vertices.
By $\B_n^{(u)}$ we denote the graded $n$-part of $\B^{(u)}$.
Of special interest is the subspace $\B^c$ (resp. $\B_n^c$ or $\B_n^{c,(u)}$)
spanned by all diagrams in $\B$ (resp. $\B_n$ or $\B_n^{(u)}$) that are
connected.

It is well-known that the primitive space in $\A$ is isomorphic to $\B^c$ and
therefore the aim of this paper is to give a lower bound for the
dimension of $\B^c_n$.

\medskip

We choose the following notation for diagrams in $\B$:
By the \asr\ we know that two univalent vertices in a nontrivial
diagram cannot be adjacent to the same trivalent vertex.
Hence we think of a connected diagram $\Gamma \in \B$ as a cubic graph, i.e. all
vertices are trivalent, 
$G(\Gamma)$ with some edges ("legs") - with a free end - attached to the edges
of $G(\Gamma)$.

In pictures the legs will be given as a number that is posed according to
the cyclic ordering at the trivalent vertices of the legs, e.g.:

\begin{displaymath}
\parbox{3cm}{\epsfbox{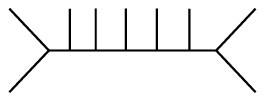}} \mapsto \parbox{3cm}
{\epsfbox{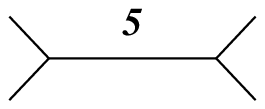}}
\end{displaymath}

\section{Weight systems coming from Lie algebras}

Let $\Lie$ be a finite dimensional Lie algebra equipped with a metric $t$,
i.e. an {\it ad}-invariant, non-degenerated, symmetric bilinear form.
There is a well-known and often described way (e.g. \cite{Kontsevich},
\cite{BarNatan1},\cite{CD2}, \cite{Vogel2} and \cite{Vaintrob2}) to use $\Lie$ for the
construction of weight systems.

A diagram $\Gamma \in \A$ will be mapped to an element $\W_{\Lie,t}(\Gamma)$
in the center $Z(U(\Lie))$ of the universal enveloping algebra of $\Lie$,
which is weighted by the word length of the elements.
$\W_{\Lie,t}(\Gamma)$ is of weight less or equal to the number of vertices
of $\Gamma$ lying on the oriented circle. This map yields a map
$\B \rightarrow Z(U(\Lie))$ (also denoted by $\W_{\Lie,t}$)
by the isomorphism $\B \rightarrow \A$.

Now let $\Lie$ be the Lie algebra $\gl(N)$, let
$e_{ij}$ be the standard generators and let
$t$ be the trace of the product of matrices.

The elements
\begin{displaymath}
c_j := \sum _{i_1, \dots, i_j} ^N e_{i_1 i_2} e_{i_2 i_3} \cdots e_{i_j i_1},
\qquad 1 \leq j \leq N,
\end{displaymath}
in $U(\gl(N))$ are called generalized Casimir elements. It is well known
(e.g. \cite {Zhelobenko}) that $Z(U(\gl(N)))$ is a free commutative polynomial
algebra in the $c_j$, $j = 1, \dots, N$.
We regard $N$ as a variable and set $c_0:=N$.

As in \cite{CD2} we will make only use of the
part $\oWgl$ of $\W_{\gl(N)}$ with highest weight; that means
for a diagram $\Gamma \in \B$ with $u$ univalent vertices that
$\oWgl(\Gamma)$ is the part of $\W_{\gl(N)}(\Gamma)$
with weight $u$. Here the weight of $c_j$ is $j$.
For a diagram $\Gamma \in \B^{(u)}$ and an arbitrary $N \geq u$ we have a nice
combinatorial description of $\oWgl$ (see \cite{CD2}
and compare with \cite{BarNatan1}):

\begin{definition}
A $\B$-state $s$ of a diagram $\Gamma \in \B$ is a map from the
internal (i.e. trivalent) vertices
of $\Gamma$ to $\{-1,1\}$. The number $\vert s \vert$ is the number of $-1$'s in $s$.

It is convenient to distinguish between two types of internal vertices:
Proper internal vertices that are not adjacent to univalent vertices
and non-proper internal vertices that are. We will denote the part
of $s$ corresponding to proper internal (resp. non-proper internal) vertices
by $s_{p}$ (resp. $s_{np}$).

Let $F(\Gamma,s)$ be the orientable surface with some missing points on the
boundary that we get by the construction:

\begin{enumerate}
\item[-] each edge of $\Gamma$ will be thickened:

\parbox{4cm}{\epsfxsize=2.5cm \epsfbox{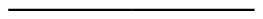}} ~
\parbox{2cm}{\epsfxsize=0.8cm \epsfbox{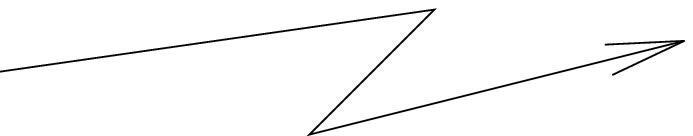}} ~
\parbox{4cm}{\epsfxsize=2.5cm \epsfbox{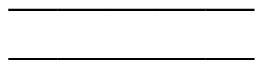}}

\item[-] each trivalent vertex will be resolved according to the value of
$s$ at it:

\begin{tabular}{lcll}
\parbox{4cm}{\epsfxsize=1.7cm \epsfbox{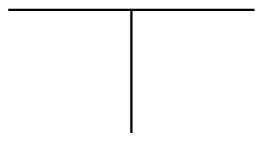}}
&\parbox{2cm}{\epsfxsize=0.8cm \epsfbox{toto.eps}} &
\parbox{4cm} {\epsfxsize=1.7cm \epsfbox{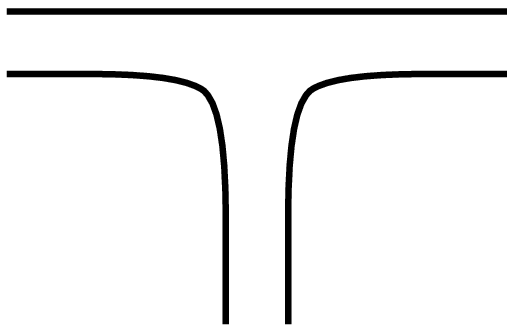}} & for $s=1$\\
&&&\\
& \parbox{2cm}{\epsfxsize=0.8cm \epsfbox{toto.eps}} &
\parbox{4cm} {\epsfxsize=1.7cm \epsfbox{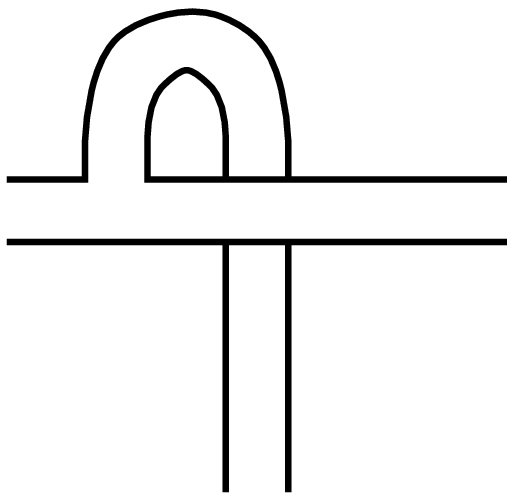}} & for $s=-1$
\end{tabular}

\item[-] each univalent vertex is responsible for a missing point in the
boundary.
\end{enumerate}

Now, for a diagram $\Gamma \in \B^{c,(u)}_n$ and a $\B$-state $s=(s_{p}, s_{np})$
we map via a function $\Omega$ the orientable surface $F(\Gamma,s)$,
to a monomial in the polynomial ring $\C[c_0, c_1, \dots, c_N]$, $N \geq u$,
in the generalized Casimir elements:
If $F(\Gamma,s)$ has
boundary components $K_1, \dots K_j$ and the number of missing points on
it are $r_1, \dots, r_j$ then $\Omega (F(\Gamma,s))$ will be
the monomial $c_{r_1} \cdots c_{r_j}$.
\end{definition}

\begin{proposition}[see \cite {CD2} and also \cite{BarNatan1}]
Let $\Gamma$ be a connected diagram in $\B^{c,(u)}_n$ - i.e.
$\Gamma$ has $(2n-u)$ trivalent vertices - and let
$N \geq u$. Then:
\begin{displaymath}
\oWgl (\Gamma) =
\sum _{{s_{p} \in \{\pm 1\} ^{2n-2 u}}}
\sum _{{s_{np} \in \{\pm 1\} ^{u}} \atop {s=(s_{p}, s_{np})}}
(-1)^{\vert s \vert}
\Omega (F(\Gamma,s)).
\end{displaymath}
\end{proposition}

The following is quite easy to see:
\begin{lemma} \label{conjugate_state}
Let $\Gamma$ be a diagram in $\B$ and let $s$ be a $\B$-state of $\Gamma$.
The conjugate $\overline s$ of $s$ is the state that we get by multiplying
each component of $s$ by $-1$.

We have:
\begin{displaymath}
F(\Gamma, s) = F(\Gamma, {\overline s}).
\end{displaymath}
\end{lemma}

We give an example for the computation of $\oWgl$:
\begin{example}
For $u$ even let $\Gamma \in \B^{c,(u)}_{u}$ be the diagram:
\begin{center}
\parbox{4cm}{\epsfxsize=1.4cm \epsfbox{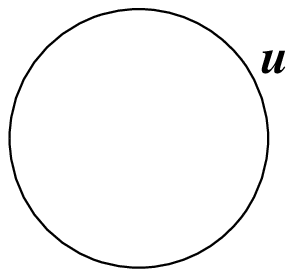}}
\end{center}
Then
$$\oWgl(\Gamma)=\sum^u_{j=0} (-1)^j {u \choose j} c_j c_{u-j}$$
\end{example}

Since in general the whole polynomial $\oWgl$  is
arduous to handle, with great success in \cite{CD2} the following
part of it is used instead:
\begin{definition}
For a diagram $\Gamma \in \B^{c,(u)}_n$, $u \leq N$, the polynomial
$CD(\Gamma)$ is the highest degree homogeneous part of
$\oWgl$, that means now the degree of each $c_j$ is one.
\end{definition}

\subsection{Evaluations for a special type of diagrams}

\begin{lemma} \label{3connected}
Let $\Gamma$ be a diagram in $\B_n^{c,(u)}$, $u$ even, $u \leq n$, 
so that the underlying cubic
graph of $\Gamma$ is planar and 3-connected. Furthermore we assume
that $\Gamma$ is
embedded in the $2$-sphere.
Then the polynomial $CD(\Gamma)$ is
\begin{eqnarray*}
CD(\Gamma) = 2 \sum_{{s_{np} \in \{\pm 1\}^{2 n-2 u}} \atop
{s_{p} = (1, \dots, 1)}, \, s=(s_p, s_{np})}
(-1)^{\vert s \vert} \Omega \circ F(\Gamma,s).
\end{eqnarray*}
\end{lemma}

\begin{proof}
(Compare with \cite{BarNatan5}.)
For a state $s$ let $\overline{F(\Gamma,s)}$ be the closed orientable
surface, obtained by gluing disks into the boundary components of
$F(\Gamma,s)$.
The Euler characteristic of $\overline{F(\Gamma,s)}$ is
\begin{displaymath}
\chi(\overline{F(\Gamma,s)})= 2 - 2 g(\overline{F(\Gamma,s)})
\end{displaymath}
where $g$ is the genus.

On the other hand we know that
\begin{displaymath}
\chi(\overline{F(\Gamma,s)}) - \# \partial(F(\Gamma,s))=
-\# e(G(\Gamma)) + \#v(G(\Gamma))
\end{displaymath}
where $\# \partial(F(\Gamma,s))$ is the number of boundary components
and $\# e(G(\Gamma))$ (resp. $\# v(G(\Gamma))$ is the number of
edges (resp. vertices) of the underlying cubic graph of $\Gamma$.

For $s=(s_p, s_{np})$ the genus $g(\overline{F(\Gamma,s)})$ only
depends on $s_p$.
Furthermore it is easy to see that
each $s_p$ that lead to a $2$-sphere $\overline{F(\Gamma,s)}$
induces an embedding of the underlying cubic graph $G(\Gamma)$ into
the $2$-sphere and vice versa.
Hence, the number of $s_p$ such that $g(\overline{F(\Gamma,s)})=0$
is equal to the number of embeddings of the underlying cubic graph
$G(\Gamma)$ into the oriented $2$-sphere.

By the classical result of Whitney we know that a $3$-connected planar
graph has only one embedding into the $2$-sphere up to homeomorphisms.
Therefore the number of $s_p$ so that $\overline{F(\Gamma,(s_p,.))}$
is the $2$-sphere is two,
i.e. corresponds to an embedding and to its mirror image.
Because the diagram is already embedded this means $s_p=(1,\dots,1)$ or
$s_p=(-1,\dots,-1)$. Using Lemma \ref{conjugate_state} we get the
desired formula.
\end{proof}

\begin{example} \label{Basic_example}
For $u=a_1+\dots+a_k+2b$ even and $n=k+u$ let
$\Gamma(a_1,\dots,a_k, b)$ be the following 'Pont Neuf' diagram in
$\B^{c,(u)}_{n}$:
\begin{center}
\parbox{3cm}{\epsfbox{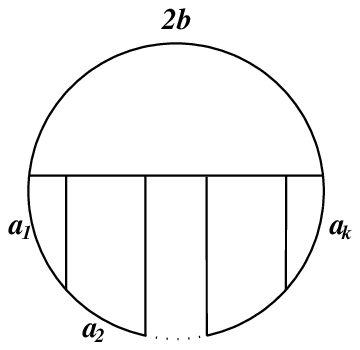}}
\end{center}
This diagram is $3$-connected and we have:
\begin{displaymath}
CD(\Gamma(a_1,\dots,a_k,b))=2 \sum_{j_1,\dots,j_k,l}^{a_1,\dots,a_k,2 b}
(-1)^{j_1+\dots j_k+l} {a_1 \choose j_1} \cdots {a_k \choose j_k}
{2 b \choose l} c_{j_1} \cdots c_{j_k} c_l c_{u-j_1-\dots-j_k-l}.
\end{displaymath}
\end{example}

\section {Asymptotic behavior of partition numbers}

Let $p(n)$ be the number of partitions of $n=a_1+ \dots + a_k$ into numbers
$1 \leq a_1 \leq a_2 \leq \dots \leq a_k$ and $p_{2}(n)$ be the number of
partitions so that $a_1\geq 2$, that means $p_2(n)=p(n)-p(n-1)$.

A theorem of Hardy and Ramanujan (see for example \cite {Hardy}), which is as
beautiful as it is famous, gives us an
asymptotic formula for $p(n)$. The asymptotic for $p_{2}(n)$ follows
by a straight forward computation and should be well-known:

\begin{theorem}[Hardy and Ramanujan] \label{HardyRamanujan}
\begin{eqnarray}
p(n) &\approx& \frac 1 {4 n \sqrt{3}} e^{\pi \sqrt{2/3} \sqrt n}\\
p_{2}(n) &\approx& \frac {\pi \sqrt{\frac 2 3}}{8 n \sqrt 3 \sqrt n}
e^{\pi \sqrt{2/ 3} \sqrt n} \, = \, \frac {\pi \sqrt{2}}{24 n \sqrt n}
e^{\pi \sqrt{2/ 3} \sqrt n}
\end{eqnarray}
\end{theorem}

In the course of this text we will need the following lemma:
\begin{lemma} \label{admissible_partitions}
Let $adm_{2}(n)$ be the number of partitions of $n=a_1 + \dots +a_k$,
so that $2 \leq a_1 \leq \dots \leq a_k$ and $n-a_k$ is even.
Then
\begin{displaymath}
\frac 1 2 \, p_{2}(n) \leq adm_{2}(n) \leq p_{2}(n).
\end{displaymath}
\end{lemma}
\begin{proof}{\bf (S. Chmutov, O.D.)}
Call a partition $(a_1, \dots, a_k)$ with $a_1 \geq 2$
admissible if $n-a_k$ is even. It is non-admissible if $n-a_k$ is odd.

Now let a partition $(a_1, \dots, a_k)$ be non-admissible.
Therefore one of $a_1, \dots, a_{k-1}$ must be odd and hence there
is an $l<k$ so that $a_l \geq 3$ and either $l=1$ or $a_{l-1}=2$.

With this choice for $l$ the map
\begin{eqnarray*}
(a_1,\dots, a_l, \dots, a_k) & \mapsto & (a_1, \dots, a_l-1, \dots , a_k+1)
\end{eqnarray*}
extends to an injective map from the set of non-admissible partitions of
cardinality $p_2(n) - adm_2(n)$ to the set of admissible ones.
Hence
\begin{displaymath}
p_2(n)- adm_2(n) \leq adm_2(n)
\end{displaymath}
\end{proof}

\section{A lower bound} \label{sectionBound}

\begin{theorem} \label{dimension_1}
For fixed $u$ and $k$ with $u$ even let ${\cal S}_{k,u}$
be the set of all $(a_1, \dots,a_k,b)$
so that
\begin{enumerate}
\item $0 \leq a_1 \leq \dots \leq a_k \leq b$
\item $u:=a_1 + \dots + a_k +2 b$.
\end {enumerate}

Let $\Gamma(a_1, \dots , a_k, b) \in \B_{u+k}^{c,(u)}$ as in Example
\ref{Basic_example}.
Then the polynomials $CD(\Gamma(a_1, \dots a_k, b))$
are linearly independent on ${\cal S}_{k,u}$.
\end{theorem}
\begin{proof}
In ${\cal S}_{k,u}$ we have an ordering by
the lexicographical ordering on $(a_1, \dots, a_k)$.

Let $(a_1,\dots,a_k,b)$ be an element of ${\cal S}_{k,u}$.
Then in $CD(\Gamma(a_1,\dots,a_k,b))$ the monomial $c_{a_1} \cdots c_{a_k}
c_b^2$ has a nontrivial coefficient and it does not occur
in any $CD(\Gamma(\tilde a_1,\dots,\tilde a_k,\tilde b))$ for
$(\tilde a_1,\dots,\tilde a_k,\tilde b)$ in ${\cal S}_{k,u}$
less than $(a_1,\dots,a_k,b)$.
\end{proof}

\begin{theorem}   \label{low_boundI}
The dimension of the subspace $\B_n^c$ generated by connected graphs in
$\B_n$ is greater than or equal to the number of partitions
\begin{eqnarray}
p_1+\dots + p_r &=& n+2 \nonumber \\
2 \leq p_1 \leq \dots &\leq& p_r \label{partition_theorem}\\
n-p_r&& \mbox{ even }. \nonumber
\end{eqnarray}
\end{theorem}

\begin{proof}
For all $u_1 \neq u_2$ it holds $\B_n^{c,(u_1)} \cap \B_n^{c,(u_2)}= \{0 \}$.
Hence, by Theorem \ref{dimension_1} we know that the dimension
of $\B^c_n$ is greater than or equal to the number of $(a_1, \dots, a_k,b)$ with
\begin{eqnarray*}
0 \leq a_1 \leq \dots \leq a_k &\leq& b\\
a_1+\dots + a_k + b +b &=&u, \quad u \mbox{ even}\\
k+u&=&n.
\end{eqnarray*}
By increasing each $a_j, j=1, \dots,k$ and $b$ by 1 we
see that the dimension of $\B_n^c$ is greater than or equal to
the number of $(a_1, \dots, a_k,b)$ with
\begin{eqnarray}
1 \leq a_1 \leq \dots \leq a_k &\leq& b \nonumber \\
a_1+\dots + a_k + b+b &=&n+2 \label{partition_set}\\
n-k && \mbox{even}. \nonumber
\end{eqnarray}

Now we will look at the Young diagram (also called Ferrers diagram)
of a partition. To each partition corresponds its conjugate partition
defined by a reflection of the diagram (see Figure \ref{fig_partition}).
\begin{figure} \label{fig_partition}
\begin{center}
\parbox{4cm} {\epsfxsize=1.6cm \epsfbox{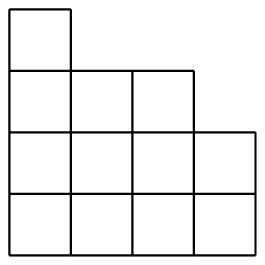}} ~~
\parbox{4cm} {\epsfxsize=1.6cm \epsfbox{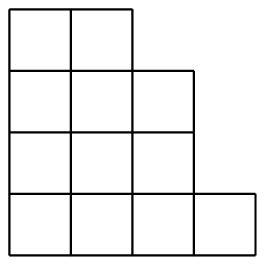}}
\end{center}
\caption{The Young diagram of the partition $12=1+3+4+4$ and
its conjugate $2+3+3+4$}
\end{figure}
By looking at the conjugates of the partitions in (\ref{partition_set})
we get the set defined by (\ref{partition_theorem}). Hence, the claim
follows.
\end{proof}

Combining Theorem \ref{low_boundI} with the Hardy-Ramanujan formula
\ref {HardyRamanujan} and Lemma \ref{admissible_partitions} we get
\begin{maintheorem}
The dimension of the primitive space in the space of Vassiliev invariants
grows in degree $n$ faster than $e^{c \sqrt n}$ for any
$c<\pi \sqrt {2/3}$, as $n$ tends to infinity.
\end{maintheorem}

\section{The dimensions of $\B_{u+k}^{c,(u)}$ for low $k$.}

For $k=0,1$ and $2$ the dimensions of $\B_{u+k}^{c,(u)}$ are known.
For $k \leq 5$ and $u$ odd we have shown in \cite{Dasbach1} that the spaces are trivial.
So for the rest of this section we assume that $u$ is even.

It is easy to see that $\B_{u}^{c,(u)}$ is one-dimensional.
Furthermore it is proved in \cite{Dasbach1} (see also \cite{Dasbach3})
that $\B_{u+1}^{c,(u)}$
is $\lfloor {\frac u 6} \rfloor +1$ dimensional.

For $k=2$ we know \cite{Dasbach3} that
$$\dim \B_{u+2}^{c,(u)} = \left \lfloor \frac {u^2+12 u} {48} \right \rfloor +1.$$

We will look at these dimension formulas from the settings given in this
paper:

\begin {enumerate}
\item {\bf $k=1:\,$} We have seen in section \ref{sectionBound} that the
dimension of $\B_{u+1}^{c,(u)}$
is greater than or equal to the number of partitions
\begin{eqnarray} \label {partition_k=1}
u & = & a_1 + 2 b, \quad 0 \leq a_1 \leq b,  \quad
a_1  \mbox{ even.}
\end{eqnarray}

With $r_1:=a_1/2$ and $r_2:=(b-a_1)$ we see that the number of partitions
fulfilling (\ref {partition_k=1}) is equal to the number of
partitions
$u = 6 r_1 + 2 r_2,\,  r_1, r_2 \geq 0$.

Hence, the generating function for the lower bound is
$$\frac 1 {(1-x^2) (1- x^6)} =  \sum_{u \,\mbox{ \tiny even}}
\left (\left \lfloor \frac u 6 \right \rfloor +1 \right ) x^u.$$

Therefore, the lower bound gives the exact dimensions.
\item {\bf $k=2: \,$} The
dimension of $\B_{u+2}^{c,(u)}$
is greater than or equal to the number of partitions $u = a_1 + a_2+ 2 b, \,
0 \leq a_1 \leq a_2 \leq b$.  Moreover $a_2 - a_1$ must be even.

With $r_1:=a_1,\, r_2:=\frac{a_2 - a_1} 2$ and $r_3:=b- a_2$ we see
that the the number of these
partitions is equal to the number of partitions
$u = 4 r_1 + 6 r_2 + 2 r_3, \, r_1, r_2, r_3 \geq 0$.

The generating function is
$$\frac 1 {(1-x^2) (1-x^4) (1- x^6)} = \sum_{u \,\mbox{ \tiny even}}
\left (\left \lfloor \frac {u^2+12 u} {48} \right \rfloor +1 \right ) x^u,$$
where one can see the equality as in \cite{Dasbach3}.

Hence, still for $k=2$ the estimate is sharp.
\item {\bf $k=3: \,$} The dimension of $\B_{u+3}^{c,(u)}$ is greater than or 
equal to the number of partitions
\begin{equation} \label {part_k3}
u=a_1+a_2+a_3+2 b, \qquad 0 \leq a_1 \leq a_2 \leq a_3 \leq b.
\end{equation}

Because $u$ is even we see that $a_1+a_3-a_2$ must be even.
The set of partitions in (\ref{part_k3}) divides into two subsets:
\begin{enumerate}
\item {\bf $a_1$ even: } In this case we set
$$r_1:=\frac {a_1} 2, \, r_2:=a_2-a_1, \,
r_3:= \frac {a_3-a_2} 2, \, r_4:=b-a_3$$
and hence we look at the number of partitions
\begin{eqnarray}
u&=& 10 r_1 + 4 r_2 + 6 r_3 + 2 r_4 \qquad r_1, r_2, r_3, r_4 \geq 0.
\end{eqnarray}
A generating function for this number is
$$\frac 1 {(1-x^2) (1-x^4) (1-x^6) (1-x^{10})}.$$
\item {\bf $a_1$ odd: } We set $\tilde a_1:=a_1-1, \, \tilde a_2:=a_2 -1,
\, \tilde a_3:=a_3-2,
\, \tilde b:=b-2.$
This yields a partition of $u-8=\tilde a_1+\tilde a_2+\tilde a_3 +2 \tilde b$
satisfying that $\tilde a_1$ is even.
Hence we have a generating function for their numbers:
$$\frac {x^8} {(1-x^2) (1-x^4) (1-x^6) (1-x^{10})}.$$
\end{enumerate}
We only have to add these two generating functions and we get a
generating function for the lower bound of the dimension of
$\B_{u+3}^{c,(u)}$ coming from our construction:
$$\frac {1+ x^8} {(1-x^2) (1-x^4) (1-x^6) (1-x^{10})}.$$
\end{enumerate}

\begin{remark}
Dror Bar-Natan \cite{BarNatan4} has computed - with the help of weight systems
coming from $so(N)$ - the dimensions of $\B^{c,(u)}_{u+3}$, $u=2,4,6,8$ and gave lower
bounds for $u=10, 12,14$. 
With a different approach Jan Kneissler
verified these dimensions \cite{Kneissler}, showed that the lower bounds are actually the
dimensions  and  
in addition gave the dimension for $u=16$:
\begin{equation} \label{BarNatanEstimates}
2,\, 3, \, 5, \, 8, \, 10, \, 15, \, 19, \, 24.
\end{equation}

Our estimate only gives $\dim \B^{c,(2)}_5 \geq 1$ and is therefore not sharp
even for the most simple case. (In fact, the reason for this is that
$CD(\Gamma(a_1,a_2,a_3,b))$ is invariant under any permutation of
$a_1, a_2$ and $a_3$, while in general the diagrams $\Gamma(a_1,a_2,a_3,b)$
are not.)

However, the estimates (\ref{BarNatanEstimates}) are the first terms of the
generating function in the 

\begin{conjecture}
The dimension of the space $\B^{c,(u)}_{u+3}$ is given by the 
generating function:
$$ \frac{1+x^2+x^8-x^{10}} {(1-x^2)(1-x^4)(1-x^6)(1-x^{10})}.$$
\end{conjecture}
\end{remark}

\section{Upper bound for the dimension in $Z(U(\gl(N)))$}
We have proved that the image of the map
$$\W_{\gl(N)}: \B_n \rightarrow Z(U(\gl(N)))$$ is greater than
$g_1(n) \, e^{c \sqrt n}$ for some rational function $g_1(n)$ 
if $N>n$ and $c = \pi \sqrt{2/3}$.

Now let $\B^r_n$ (resp. $\B^{r,(u)}_n$) be the subspace of $\B_n$
(resp. $\B^{(u)}_n$) that is generated by diagrams with at least one
trivalent vertex in each connected component.
Let $\Gamma$ be a diagram in $\B^{r,(u)}_n$. Specifically this means
that $u \leq n$.

$\W_{\gl(N)}(\Gamma)$, $N>u$, is a polynomial in the
generalized Casimir elements $c_1, \dots, c_u$ and $c_0:=N$.

A monomial $c_{j_1} \cdots c_{j_r}$ in $\W_{\gl(N)}(\Gamma)$
fulfills $\sum _{k=1}^{r} j_k \leq u$. Furthermore as in the proof of
Lemma \ref{3connected} we know that the homogenous degree of the polynomial
$\W_{\gl(N)}(\Gamma)$ is less or equal to $n$.

Therefore, $\dim \W_{\gl(N)} (\B^r_n)$ is less or equal to the number of
partitions
\begin{eqnarray*}
0 \leq j_1 \leq \dots & \leq& j_r\\
j_1+\dots+j_r & \leq & n\\
r&\leq & n.
\end{eqnarray*}

A rough estimate gives:
\begin{lemma}
$$\dim \W_{\gl(N)}(\B^r_n) \leq n^2 p(n),$$
where $p(n)$ is the partition number.
\end{lemma}

\begin{remark}
In fact, we only used the space $\B^r_n$ instead of $\B_n$ itself
to avoid some messy details. Let $l$ be the element in $\B$ with two
univalent vertices and no trivalent vertex. $l^k$ is the disjoint union of
$k$ copies of $l$.

So $\B_n = \bigoplus_{k=0}^n l^k \B_{n-k}^r$.
Looking carefully at the isomorphism of $\B$ and $\A$ one can
see that
$$\dim \W_{\gl(N)}(\B_n) \leq \sum_{k=0}^n \dim \W_{\gl(N)}(\B^r_{n-k})
\leq n^3 p(n).$$
\end{remark}

Recently there was some interest in  operations on the space of Vassiliev invariants
that where induced by doing cabling operations
(see e.g. \cite{KSA}, \cite{McDR} or \cite{CDV}).
One can interpret the facts given in this section as (see \cite{CDV} for
details):
\begin{corollary}
Let $C_a H(n)$ be the space of all Vassiliev invariants coming from the
HOMFLY-Vassiliev invariants and all of their cablings, connected and
disconnected.
Then there are rational functions $g_1(n)$ and $g_2(n)$ so that
$$g_1(n) e^{c \sqrt n} \leq \dim C_a H(n) \leq g_2(n) e^{c \sqrt n}$$.

\end{corollary}

\providecommand{\bysame}{\leavevmode\hbox to3em{\hrulefill}\thinspace}

\end{document}